\documentclass[12pt]{amsart}
\usepackage{a4wide}
\usepackage[margin=0.8in]{geometry}
\usepackage{amsmath,amssymb}
\usepackage{tikz}
\usepackage{tikz-cd}
\usetikzlibrary{arrows.meta}

\usepackage{amsthm}
\usepackage{mathtools}

\usepackage{color}
\usepackage{upref}
\usepackage{enumitem}
\usepackage{verbatim}
\usepackage{mathrsfs}
\usepackage{hyperref}

\allowdisplaybreaks

\newcommand{\AXA}{$\pre AX_A$}

\newcommand{\ASB}{$\pre AS_B$}

\newcommand{\N}{\mathbb{N}}

\DeclareMathOperator{\spn}{span}

\newcommand{\ip}[2]{\langle #1,#2\rangle}

\allowdisplaybreaks

\newcommand{\<}{\langle}
\renewcommand{\>}{\rangle}

\theoremstyle{plain}
\newtheorem{theorem}{Theorem}[section]
\newtheorem{proposition}[theorem]{Proposition}
\newtheorem{lemma}[theorem]{Lemma}

\theoremstyle{definition}

\newtheorem{remark}[theorem]{Remark}

\numberwithin{equation}{section}

\numberwithin{equation}{section}

\newcommand{\KK}{\mathcal K}

\newcommand{\LL}{\mathcal L}
\newcommand{\pre}[1]{{}_{#1}}

\setlist{font=\normalfont} 
\newcommand{\ots}{\otimes_A}
\newcommand{\otss}{\otimes_B}
\newcommand{\tpsi}{\tilde{\psi}}
\newcommand{\tvarphi}{\tilde{\varphi}}
\newcommand{\trho}{\tilde{\rho}}

\newcommand{\norm}[1]{\left\lVert#1\right\rVert}

\begin{document}

\title{EQUIVALENCES OF $\textrm{C}^*$-CORRESPONDENCE RELATIONS}

\author[Boyu Li]{Boyu Li}
\address{Department of Mathematical Sciences
\\New Mexico State University
\\Las Cruces, New Mexico 88003}
\email{boyuli@nmsu.edu}

\author{Menev\c{s}e  Ery\"uzl\"u Paulovicks}
\address{Department of Mathematical Sciences
\\New Mexico State University
\\Las Cruces, New Mexico 88003}
\email{meryuzlu@nmsu.edu}



\date{\today}

\thanks{B. Li is partially supported by an NSF grant (DMS-2350543).}

\begin{abstract}
We study three relations between $\textrm{C}^*$-correspondences that lie between
strong shift equivalence and shift equivalence, called aligned, balanced, and compatible shift equivalence. We show that for regular and full $\textrm{C}^*$-correspondences, the three relations
coincide. 
\end{abstract}

\maketitle

\section{Introduction}\label{sec:intro}

In his foundational work,  Williams characterized topological conjugacy of shifts of finite type in terms of strong shift equivalence of their adjacency matrices and introduced shift equivalence as a potentially more tractable algebraic relation \cite{Wil73}. Two square matrices $A$ and $B$ with entries in $\mathbb{N}$ are called shift equivalent with lag $m$ if there are rectangular
matrices $R$ and $S$ with entries in $\mathbb{N}$ such that
\[
A^m = RS, \qquad B^m = SR, \qquad AR = RB, \qquad SA = BS.
\]
They are called strong shift equivalent if they are related by a finite chain of
shift equivalences with lag one. The distinction between shift equivalence and strong shift equivalence subsequently became a central problem in the study of subshifts, and Williams conjectured that the two relations 
coincide. However, this was disproved by Kim and Roush \cite{KR92} for reducible shifts and \cite{KR99} for irreducible shifts. 

In 2008, Muhly, Pask, and Tomforde \cite{MPT08} extended Williams' program to $\textrm{C}^*$-correspondences, motivated by the graph correspondence studied by Cuntz and Krieger \cite{CK80}. It was observed in \cite{CK80} that the Cuntz-Krieger $\textrm{C}^*$-algebras of strong shift equivalent matrices are Morita equivalent. Building upon this observation, one can extend William's definitions of shift and strong shift equivalences of matrices to $\textrm{C}^*$-correspondences, where matrix multiplications are replaced by balanced tensor products. It is proved in \cite{MPT08} that, similar to the matrix case, the Cuntz-Pimsner $\textrm{C}^*$-algebras of strong shift equivalent $\textrm{C}^*$-correspondences are Morita equivalent. As a result, it is natural to consider the implications of shift and strong shift equivalence in this operator-algebraic setting, in particular, on the Cuntz-Pimsner algebra associated with $\textrm{C}^*$-correspondences \cite{AHL2022, BDR-hom, BDR-SER, BMR24, CDE24, Ery22, EKK2017, KK14, KK2022corr,  MS2019}. In particular, whether shift equivalent $\textrm{C}^*$-correspondences have Morita equivalent Cuntz-Pimsner algebras remains an open problem. 

To clarify the relationship between shift and strong shift equivalences, Carlsen, Dor-On, and Eilers introduced compatible shift equivalence \cite{CDE24}. For $\mathbb{N}$-valued essential matrices, compatible shift equivalence is equivalent to strong shift equivalence. Subsequently, Bilich, Dor-On, and Ruiz \cite{BDR-SER} introduced aligned and balanced shift equivalences of matrices, which are also equivalent to strong shift equivalence. As aligned, balanced, and compatible shift equivalences are all equivalent to strong shift equivalence, they provide alternative ways to computationally verify whether two matrices are strong shift equivalent.

The notions of aligned, balanced, and compatible shift equivalences can easily be extended to $\textrm{C}^*$-correspondences. 
However, it is unclear how they are related for general $\textrm{C}^*$-correspondences. 
The proof that aligned, balanced, and compatible shift equivalences of matrices are all equivalent in \cite{BDR-SER} relies on combinatorial arguments involving bijections between sets of paths, which do not directly extend to general $\textrm{C}^*$-correspondences. 

In this paper, we prove that aligned, balanced, and compatible shift equivalences among regular and full $\textrm{C}^*$-correspondences are equivalent (Theorem~\ref{main}). One of the key techniques is Proposition~\ref{prop.alpha.map}, which allows us to recover a correspondence automorphism from compatible factorizations on tensor products of $\textrm{C}^*$-correspondences.

We would like to thank Adam Dor-On for many discussions on this topic and for sketching a proof for Proposition~\ref{prop.alpha.map}. 


\section{Preliminaries}

A $\textrm{C}^*$-correspondence $\pre AX_B$ is a right Hilbert $B$-module equipped with a left action given by a homomorphism $\varphi_X: A\rightarrow \LL(X),$ where $\LL(X)$ denotes the $\textrm{C}^*$-algebra of adjointable operators on $X$. The correspondence $\pre AX_B$ is called \emph{nondegenerate} when $A\cdot X=X$. We assume throughout this paper that \emph{all $\textrm{C}^*$-correspondences are nondegenerate}.

A $\textrm{C}^*$-correspondence $\pre AX_B$ is called \emph{injective} if the left action $\varphi_X: A\rightarrow \LL(X)$ is injective; it is called \emph{proper} if  $\varphi_{X}(A)$ is contained in the $\textrm{C}^*$-algebra $\KK(X)$ of compact operators on $X$. A $\textrm{C}^*$-correspondence $\pre AX_B$ is called \emph{regular} if it is both injective and proper. It is called \emph{full} if $\overline{\text{span}}\{\<x,y\>_B: x,y\in X\}=B$

A $\textrm{C}^*$-correspondence isomorphism is a triple $(\Phi, \iota_A, \iota_B): \pre AX_B \rightarrow \pre AY_B$ consisting of a bijective linear map $\Phi: X\rightarrow Y$, and $*$-isomorphisms $\iota_A: A\rightarrow A$ and $\iota_B: B\rightarrow B$   satisfying
\begin{enumerate}
\item $\Phi(a\cdot x)=\iota_A(a)\cdot \Phi(x)$,
\item $\iota_B( \<x,z\>_B ) = \<\Phi(x), \Phi(z)\>_{B}$,
\end{enumerate}
for all $a\in A$, and $x,z\in X.$ In this paper, we take $\iota_A$ and $\iota_B$ to be the identity maps, and denote the $\textrm{C}^*$-correspondence isomorphism simply by $\Phi$.

The \emph{balanced tensor product} $X\otimes_BY$ of
an $A-B$ correspondence $X$ and a $B-C$ correspondence $Y$ is
formed as follows:
the algebraic tensor product $X\odot Y$
has an $A-C$ bimodule structure satisfying
\[
a(x\otimes y)c=ax\otimes yc,
\] 
for $a\in A,x\in X,y\in Y,c\in C$, and the unique $C$-valued semi-inner product whose values on elementary tensors are given by
\[
\<x\otimes y,u\otimes v\>_C=\<y,\<x,u\>_B\cdot v\>_C\]
for $x,u\in X,$ and $y,v\in Y$.

This semi-inner product defines a $C$-valued inner product on the quotient $X{\odot}_BY$ of $X\odot Y$ by the subspace generated by elements of form 
\[ x\cdot b \otimes y - x\otimes \varphi_Y(b)y, \]
for $b\in B$. The completion $X\otimes_B Y$ of  $X{\odot}_BY$ with respect to the norm coming from the $C$-valued inner product is an $A-C$ correspondence, where the left action is given by 
\[A\rightarrow \LL(X\otss Y), \hspace{.5cm} a\mapsto \varphi_X(a)\otimes 1_Y.\] 
We denote the canonical image of $x\otimes y$ in $X\otss Y$ by $x\otss y$. The term \emph{balanced} refers to the property
\[
x\cdot b\otss y=x\otss b\cdot y
\hspace{.5cm} \text{for } x\in X,b\in B,y\in Y,
\]
which is a consequence of the construction.

\section{Technical Lemmas}

We start this section with the following two technical lemmas. The first one follows immediately from the injectivity of the left $A$-action. 

\begin{lemma}\label{essential}
Let \AXA , \ASB\ be $\textrm{C}^*$-correspondences and let $T\in \LL(X).$ If S is injective, then the map $T\mapsto T\otimes 1_S$ is injective.
\end{lemma}

\begin{lemma}\label{essential2}
Let $\pre AX_A, \pre AY_A$ be $\textrm{C}^*$-correspondences. Let $U: Y\rightarrow Y$ be a $\textrm{C}^*$-correspondence automorphism. If $\pre AX_A$ is full and $1_X\otimes U=1_X\otimes 1_Y,$ then $U$ is the identity may on $Y$.
\end{lemma}

\begin{proof}
Let $y,y'\in Y.$ Then, Cohet-Hewitt factorization theorem and fullness of $\pre AX_A$ together implies $y'=\lim_{\alpha}\sum_{i=1}^{N_n}\<x_{\alpha,i},x'_{\alpha,i}\>_A\cdot y''$, for  $x_{\alpha}, x'_{\alpha}\in X,$ and $y''\in Y.$ Then we have 
\begin{align*}
\<U(y),y'\>_A&=\lim_\alpha\sum_{i=1}^{N_n}\<U(y),\<x_{\alpha,i},x'_{\alpha,i}\>_A\cdot y''\>_A\\
&=\lim_{\alpha}\sum_{i=1}^{N_n}\<x_{\alpha,i}\ots U(y),x'_{\alpha,i}\ots y''\>_A\\
&=\lim_{\alpha}\sum_{i=1}^{N_n}\<x_{\alpha,i}\ots y,x'_{\alpha,i}\ots y''\>_A\\
&=\<y,\lim_{\alpha}\sum_{i=1}^{N_n}\<x_{\alpha,i},x'_{\alpha,i}\>_A\cdot y''\>_A\\
&=\<y,y'\>_A,
\end{align*}
completing the proof.\end{proof}

We now state a key technique of this paper.

\begin{proposition}\label{prop.alpha.map}
Suppose ${}_AX_A$, ${}_BY_B$, and
${}_AS_B$ are regular and full $\textrm{C}^*$-correspondences. Suppose
\mbox{$T\in\LL(X\otimes_A S\otimes_B Y)$} is a unitary such that
\[
T = 1_X\otimes\tau_1 = \tau_2\otimes 1_Y ,
\]
where $\tau_1\in\LL(S\otimes_B Y)$ is a left $A$-linear unitary and
$\tau_2\in\LL(X\otimes_A S)$ is a unitary. Then there is a unique $\textrm{C}^*$-correspondence automorphism $\tau_3$ of $\pre AS_B$ such that
$\tau_1=\tau_3\otimes 1_Y$. In particular $T = 1_X\otimes\tau_3\otimes 1_Y$.
\end{proposition}

\begin{proof}
For $s\in S$ define a linear map  \mbox{$L_s: Y\rightarrow  S\otimes_B Y$} on elementary tensors by  $L_s(y)=s\otimes_B y$. It is easy to verify that this map is adjointable with  
\[L_s^*(s'\otimes_B y)=\varphi_Y(\<s,s'\>_B)y, \text{ and } L_s^*L_{s'} = \varphi_Y(\<s,s'\>_B).\] Then we have $\norm{L_s}^2=\norm{\varphi_Y(\<s,s\>_B)}=\norm{s}^2.$ Thus the map $s\mapsto L_s$ is a linear
isometry \mbox{$S\rightarrow\LL(Y,S\otimes_B Y)$.} Moreover, we have the  equalities
\begin{equation}\label{L-equalities}
L_{sb}=L_s\varphi_Y(b), \text{ and } L_{a\cdot s}=(\varphi_S(a)\otimes 1_Y)L_s,
\end{equation}
for $a\in A,\ b\in B,\ s\in S$. Similarly, for $x\in X$ define the map $L_x: S\rightarrow X\ots S$, $L_x(s)=x\ots s$.
Notice that the map $M_x: S\otss  Y\to X\ots S\otss Y$ defined by  $M_x(w)=x\ots w$
satisfies \[M_x = L_x\otimes 1_Y \text{ and } M_x^*(x'\ots w)=(\varphi_S(\<x,x'\>_A)\otimes 1_Y)w,\] for any $x\in X, w\in S\otss Y.$ Now let $x,x'\in X$. On one hand we have 
\[M_x^*(1_X\otimes\tau_1)M_{x'}=[\varphi_S(\<x,x'\>_A)\otimes 1_Y]\tau_1, \]
by definition. On the other hand,
\[
M_x^*(\tau_2\otimes 1_Y)M_{x'}
=(L_x^*\otimes 1_Y)(\tau_2\otimes 1_Y)(L_{x'}\otimes 1_Y)
=(L_x^*\tau_2L_{x'})\otimes 1_Y .
\]
Since $1_X\otimes\tau_1=\tau_2\otimes 1_Y$, we conclude
\begin{equation}\label{key_equation}
(L_x^*\tau_2L_{x'})\otimes 1_Y=[\varphi_S(\ip{x}{x'}_A)\otimes 1_Y]\tau_1,
\end{equation}
for any $x,x'\in X.$

We are now ready to construct the automorphism $\tau_3$.
Choose a positive contractive approximate identity $(e_\lambda)_{\lambda\in\Lambda}$ for $A$. Since $X$ is full, $I=\spn\{\<x,x'\>_A:x,x'\in X\}$ is dense in $A$. For each $(\lambda,n)\in \Lambda\times\N$, with $n\geq1$, choose $a_{\lambda,n}\in I$ such that
\[\norm{a_{\lambda,n}-e_\lambda}<\frac{1}{n}.\]
Equip $D=\Lambda\times\mathbb N$ with the product order:
$(\lambda,n)\leq(\mu,m)\iff \lambda\leq\mu \text{ and } n\leq m.$ We regard these elements as a net $(a_d)_{d\in D}$. Then for any $c\in A$ we have 
\[\norm{a_{\lambda,n}c-c}=\norm{a_{\lambda,n}c-c+e_\lambda c-e_{\lambda}c}< \norm{(a_{\lambda,n}-e_\lambda)c}+\norm{e_\lambda c-c} < \frac{\norm{c}}{n} +\norm{e_\lambda c-c}\rightarrow 0.\] 
By Cohen-Hewitt factorization theorem, we may find $c\in A, s'\in S$ such that $s=c\cdot s'.$ Then, 
\[\norm{a_d\cdot s-s}=\norm{a_d\cdot c\cdot s'-c\cdot s'}\leq\norm{a_d c-c}\norm{s'}\rightarrow 0,\] which means
$a_d \cdot s\to s$, for $s\in S.$ Write $a_d=\sum_{i=1}^{N_d}\<x_{d,i},y_{d,i}\>_A$ and define
\[
\tau_{3,d}:=\sum_{i=1}^{N_d}L_{x_{d,i}}^*\,\tau_2\,L_{y_{d,i}}\in\LL(S).
\]
Using \eqref{key_equation} we get  $\tau_{3,d}\otimes 1_Y=(\varphi_S(a_d)\otimes 1_Y)\tau_1$.
Since $\tau_1$ is left $A$-linear, we have $[\varphi_S(a_d)\otimes 1_Y]\tau_1=\tau_1[\varphi_S(a_d)\otimes 1_Y]$. Hence,
for $s\in S$ and $t_d:=\tau_{3,d}(s)$, using \eqref{L-equalities},
\begin{equation}\label{L-equalities2}
L_{t_d}=(\tau_{3,d}\otimes 1_Y)L_s=\tau_1(\varphi_S(a_d)\otimes 1_Y)L_s
=\tau_1 L_{a_d\cdot s}.
\end{equation}
As $\tau_1$ is unitary and $s\mapsto L_s$ is isometric, we have 
\[\norm{t_d-t_{d'}}=\norm{\tau_1L_{(a_d-a_{d'})\cdot s}}=\norm{(a_d-a_{d'})\cdot s}\rightarrow 0,\]
for $d,d'\in D.$ Thus $(t_d)$ is Cauchy, and we define
$\tau_3(s):=\lim_d\tau_{3,d}(s)$. Then, by  \eqref{L-equalities2} we get
\begin{equation}\label{L-equalities3}
L_{\tau_3(s)}=\tau_1 L_s, \hspace{.5cm} s\in S .
\end{equation}
Since $s\mapsto L_s$ is injective, \eqref{L-equalities3} determines $\tau_3$ uniquely. It remains to show that $\tau_3$ is a $\textrm{C}^*$-correspondence isomorphism. 
\begin{enumerate}[label=(\roman*)]
\item \emph{Left $A$-linearity:} By \eqref{L-equalities} and $A$-linearity of $\tau_1$, we have 
\[L_{\tau_3(a\cdot s)}=\tau_1(\varphi_S(a)\otimes 1_Y)L_s=(\varphi_S(a)\otimes 1_Y)\tau_1L_s=L_{a\cdot\tau_3(s)},\]
for $a\in A, s\in S.$

\item \emph{Inner-products:} For $s,s'\in S$ we compute
\[\varphi_Y(\<\tau_3 (s),\tau_3 (s')\>_B)=L_{\tau_3(s)}^*L_{\tau_3(s')}=L_s^*\tau_1^*\tau_1 L_{s'}=L_s^*L_{s'}
=\varphi_Y(\<s,s'\>_B).\] Since $\varphi_Y$ is injective, we may conclude $\<\tau_3 (s),\tau_3 (s')\>_B=\<s,s'>_B$.

\item \emph{Surjectivity:} Taking adjoins in our hypothesis equalities $T=1_X\otimes\tau_1=\tau_2\otimes 1_Y$ we get 
\mbox{$T^*=1_X\otimes\tau_1^*=\tau_2^*\otimes 1_Y$.} Now going through the same construction procedure as $\tau_3$ by using $\tau_2^*$ instead of $\tau_2$ we may define a linear map $\sigma: S\rightarrow S$ such that $L_{\sigma(s)}=\tau_1^*L_s.$ Then we have 
\[ L_s=\tau_1\tau_1^*L_s=\tau_1 L_{\sigma(s)}=L_{\tau_3(\sigma(s))},\] for any $s\in S$. Surjectivity now follows from the fact that the map $s\mapsto L_s$ is injective.  
\end{enumerate}

Now using  \eqref{L-equalities3}, for any $s\in S, y\in Y$,  we get  \[(\tau_3\otimes 1_Y)(s\otss y)=\tau_3(s)\otss y= L_{\tau_3(s)}(y)
=\tau_1 L_s(y)=\tau_1(s\otss y).\] Consequently, we have the equality $\tau_3\otimes 1_Y=\tau_1,$ which implies  \mbox{$T=1_X\otimes\tau_3\otimes 1_Y$,} as desired. 
\end{proof}

\section{Equivalences of $\textrm{C}^*$-correspondence relations}

Recall, two $\textrm{C}^*$-correspondences $\pre AX_A$ and $\pre BY_B$ are shift equivalent with lag $m$ if there exists two $\textrm{C}^*$-correspondences $\pre A R_B$ and $\pre B S_A$, and unitary maps:
\begin{align*}
    \psi_X: & R\otimes_B S\to X^{\otimes m} \\
    \psi_Y: & S\otimes_A R\to Y^{\otimes m} \\
    \varphi_R: & X\otimes_A R\to R\otimes_B Y \\
    \varphi_S: & Y\otimes_B S\to S\otimes_A X
\end{align*}

Furthermore, we say this shift equivalence is
\begin{enumerate}
    \item Aligned with lag $m$: if 
    \begin{align*}
        (\psi_X^{-1}\otimes i_{X})(i_{X}\otimes\psi_X) &=(i_R\otimes \varphi_S)(\varphi_R\otimes i_S) \\
        (\psi_Y^{-1}\otimes i_{Y})(i_{Y}\otimes\psi_Y) &=(i_S\otimes \varphi_R)(\varphi_S\otimes i_R)
    \end{align*}
    \item Balanced with lag $m$ if 
    \begin{align*}
        (\psi_X^{-1}\otimes i_{X^{\otimes m}})(i_{X^{\otimes m}}\otimes\psi_X) &=(i_R\otimes \varphi_S^{(m)})(\varphi_R^{(m)}\otimes i_S) \\
        (\psi_Y^{-1}\otimes i_{Y^{\otimes m}})(i_{Y^{\otimes m}}\otimes\psi_Y) &=(i_S\otimes \varphi_R^{(m)})(\varphi_S^{(m)}\otimes i_R)
    \end{align*}
    \item Compatible with lag $m$ if 
    \begin{align*}
        \varphi_R^{(m)} &= (i_R\otimes\psi_Y)(\psi_X^{-1}\otimes i_R) \\
        \varphi_S^{(m)} &= (i_S\otimes\psi_X)(\psi_Y^{-1}\otimes i_S)
    \end{align*}
\end{enumerate}
Here, for $\ell\geq 1$, we use the following notations:
\begin{align*}
\varphi_R^{(\ell)}&=(\varphi_R\otimes i_{Y^{\otimes (\ell-1)}})(i_X\otimes \varphi_R\otimes i_{Y^{\otimes (\ell-2)}})\cdots(i_{X^{{\otimes (\ell-1)}}}\otimes \varphi_R): X^{\otimes \ell}\otimes R\to R\otimes Y^{\otimes \ell}, \\
\varphi_S^{(\ell)}&=(\varphi_S\otimes i_{X^{\otimes (\ell-1)}})(i_Y\otimes \varphi_S\otimes i_{X^{\otimes (\ell-2)}})\cdots(i_{Y^{{\otimes (\ell-1)}}}\otimes \varphi_S): Y^{\otimes \ell}\otimes S\to S\otimes X^{\otimes \ell}.
\end{align*}
We say $X$ and $Y$ are aligned, balanced, compatible shift equivalent if they are aligned, balanced, compatible shift equivalent with some lag $m$, respectively. The main goal of this section is to prove that aligned, balanced, and compatible shift equivalences are, in fact, equivalent. 

\begin{theorem}\label{main} Let $\pre AX_A, \pre BY_B$ be two regular, full $\textrm{C}^*$-correspondences. Then the following are equivalent:
\begin{enumerate}
    \item $X$ and $Y$ are aligned shift equivalent.
    \item $X$ and $Y$ are balanced shift equivalent.
    \item $X$ and $Y$ are compatible shift equivalent.
\end{enumerate}
\end{theorem}

We shall proceed as follows: the direction (1)$\implies$ (2) is Proposition~\ref{prop.a.to.b}, (2)$\implies$(3) is Proposition~\ref{prop.b.to.c}, and (3)$\implies$(1) is Proposition~\ref{prop.c.to.a}.

\begin{proposition}\label{prop.a.to.b} Aligned shift equivalence with lag $m$ implies balanced shift equivalence with lag $m$. 
\end{proposition}

\begin{proof} Assuming aligned shift equivalence with lag $m$, we prove, inductively, that for all $\ell\geq 1$, 
\begin{align*}
        (\psi_X^{-1}\otimes i_{X^{\otimes \ell}})(i_{X^{\otimes \ell}}\otimes\psi_X) &=(i_R\otimes \varphi_S^{(\ell)})(\varphi_R^{(\ell)}\otimes i_S) \\
        (\psi_Y^{-1}\otimes i_{Y^{\otimes \ell}})(i_{Y^{\otimes \ell}}\otimes\psi_Y) &=(i_S\otimes \varphi_R^{(\ell)})(\varphi_S^{(\ell)}\otimes i_R)
\end{align*}
For $\ell=1$, this is precisely the conditions for aligned shift equivalence. Now, assume the equations are true for $\ell\geq 1$, let us prove the case for $\ell+1$. 
\begin{align*}
    &(i_R\otimes \varphi_S^{(\ell+1)})(\varphi_R^{(\ell+1)}\otimes i_S)  \\
    =& (i_R\otimes \varphi_S\otimes i_{X^{\otimes \ell}}) (i_{R\otimes Y}\otimes \varphi_S^{(\ell)})(\varphi_R\otimes i_{Y^{\otimes \ell}\otimes S})(i_X\otimes \varphi_R^{(\ell)}\otimes i_S)\\
    =& (i_R\otimes \varphi_S\otimes i_{X^{\otimes \ell}}) (\varphi_R\otimes i_{S\otimes X^{\otimes \ell}})(i_{X\otimes R}\otimes \varphi_S^{(\ell)})(i_X\otimes \varphi_R^{(\ell)}\otimes i_S)\\
    =& \Big[\big((i_R\otimes\varphi_S)(\varphi_R\otimes i_S)\big)\otimes i_{X^{\otimes \ell}}\Big] \Big[i_X\otimes \big((i_R\otimes\varphi_S^{(\ell)})(\varphi_R^{(\ell)}\otimes i_S)\big)\Big] \\
    =& \Big[\big((\psi_X^{-1}\otimes i_X)(i_X\otimes \psi_X)\big)\otimes i_{X^{\otimes \ell}}\Big] \Big[i_X\otimes \big((\psi_X^{-1}\otimes i_{X^{\otimes \ell}})(i_{X^{\otimes \ell}}\otimes \psi_X)\big)\Big] \\
    =& (\psi_X^{-1}\otimes i_{X^{\otimes (\ell+1)}})(i_{X^{\otimes (\ell+1)}}\otimes \psi_X).
\end{align*}
The other equality follows similarly. Now, set $\ell=m$, these equations are precisely the conditions for balanced shift equivalence. 
\end{proof}

\begin{lemma}\label{lemma.balanced.eq} Balanced shift equivalence with lag $m$ implies:
\begin{align*}
    \left[(\psi_Y\otimes i_S)(i_S\otimes\psi_X^{-1})\varphi_S^{(m)}\right]\otimes i_{X^{\otimes m}} &= i_{Y^{\otimes m}} \otimes \left[\varphi_S^{(m)}(\psi_Y\otimes i_S)(i_S\otimes\psi_X^{-1})\right], \\
    \left[(\psi_X\otimes i_R)(i_R\otimes\psi_Y^{-1})\varphi_R^{(m)}\right]\otimes i_{Y^{\otimes m}} &= i_{X^{\otimes m}} \otimes \left[\varphi_R^{(m)}(\psi_X\otimes i_R)(i_R\otimes\psi_Y^{-1})\right].
\end{align*}
\end{lemma}

\begin{proof} We shall prove the first equality and the other one would follow similarly. Starting from the left hand side,
\begin{align*}
    & \left[(\psi_Y\otimes i_S)(i_S\otimes\psi_X^{-1})\varphi_S^{(m)}\right]\otimes i_{X^{\otimes m}}\\
    =& (\psi_Y\otimes i_{S\otimes X^{\otimes m}})\left[i_S\otimes \left((i_R\otimes\varphi_S^{(m)})(\varphi_R^{(m)}\otimes i_S)(i_{X^{\otimes m}}\otimes \psi_X^{-1})\right)\right] (\varphi_S^{(m)}\otimes i_{X^{\otimes m}}) \\
    =& (\psi_Y\otimes i_{S\otimes X^{\otimes m}})(i_{S\otimes R}\otimes\varphi_S^{(m)})(i_S\otimes \varphi_R^{(m)}\otimes i_S)(\varphi_S^{(m)}\otimes i_{R\otimes S}) (i_{Y^{\otimes m}\otimes S}\otimes \psi_X^{-1}) \\
    =& (i_{Y^{\otimes m}}\otimes\varphi_S^{(m)})
    \left[\left((\psi_Y\otimes i_{Y^{\otimes m}})(i_S\otimes\varphi_R^{(m)})(\varphi_S^{(m)}\otimes i_R)\right)\otimes i_S\right] (i_{Y^{\otimes m}\otimes S}\otimes \psi_X^{-1}) \\
    =& (i_{Y^{\otimes m}}\otimes\varphi_S^{(m)})
    (i_{Y^{\otimes m}}\otimes \psi_Y\otimes i_S) (i_{Y^{\otimes m}\otimes S}\otimes \psi_X^{-1}) \\
    =& i_{Y^{\otimes m}} \otimes \left[\varphi_S^{(m)}(\psi_Y\otimes i_S)(i_S\otimes\psi_X^{-1})\right]. \qedhere
\end{align*}
    
\end{proof}

\begin{proposition}\label{prop.b.to.c} Balanced shift equivalence with lag $m$ implies compatible shift equivalence with lag $m$. 
\end{proposition}

\begin{proof} Apply Proposition~\ref{prop.alpha.map} to the equalities obtained in Lemma~\ref{lemma.balanced.eq}, we find two unitary maps $\alpha_S:S\to S$ and $\alpha_R:R\to R$ such that
\begin{align*}
    \left[(\psi_Y\otimes i_S)(i_S\otimes\psi_X^{-1})\varphi_S^{(m)}\right]\otimes i_{X^{\otimes m}} &= i_{Y^{\otimes m}} \otimes \left[\varphi_S^{(m)}(\psi_Y\otimes i_S)(i_S\otimes\psi_X^{-1})\right]\\
    &=i_{Y^{\otimes m}}\otimes \alpha_S\otimes i_{X^{\otimes m}}, \\
    \left[(\psi_X\otimes i_R)(i_R\otimes\psi_Y^{-1})\varphi_R^{(m)}\right]\otimes i_{Y^{\otimes m}} &= i_{X^{\otimes m}} \otimes \left[\varphi_R^{(m)}(\psi_X\otimes i_R)(i_R\otimes\psi_Y^{-1})\right]\\
    &=i_{X^{\otimes m}}\otimes \alpha_R\otimes i_{Y^{\otimes m}}.
\end{align*}
We have that, by Lemma~\ref{essential},
\begin{align}
    \label{eq.2.1}
    (\psi_Y\otimes i_S)(i_S\otimes\psi_X^{-1})\varphi_S^{(m)} &= i_{Y^{\otimes m}} \otimes \alpha_S, \\
    \label{eq.2.2}
    \varphi_S^{(m)}(\psi_Y\otimes i_S)(i_S\otimes\psi_X^{-1}) &= \alpha_S\otimes i_{X^{\otimes m}}.
\end{align}
Similarly,
\begin{align}
    \label{eq.2.3}
    (\psi_X\otimes i_R)(i_R\otimes\psi_Y^{-1})\varphi_R^{(m)} &= i_{X^{\otimes m}} \otimes \alpha_R, \\
    \label{eq.2.4}
    \varphi_R^{(m)}(\psi_X\otimes i_R)(i_R\otimes\psi_Y^{-1}) &= \alpha_R\otimes i_{Y^{\otimes m}}.
\end{align}
From Equations~\eqref{eq.2.1} and~\eqref{eq.2.3}, we have
\begin{align*}
    \varphi_S^{(m)} &= (i_S\otimes \psi_X)(\psi_Y^{-1}\otimes \alpha_S) \\
    \varphi_R^{(m)} &= (i_R\otimes \psi_Y)(\psi_X^{-1}\otimes \alpha_R) 
\end{align*}
Plugin these into the balanced shift equivalence condition $\psi_X^{-1}\otimes\psi_X=(i_R\otimes \varphi_S^{(m)})(\varphi_R^{(m)}\otimes i_S)$, we get
\begin{align*}
    \psi_X^{-1}\otimes\psi_X&=(i_R\otimes \varphi_S^{(m)})(\varphi_R^{(m)}\otimes i_S)\\ &= (i_{R\otimes S}\otimes \psi_X)(i_R\otimes \psi_Y^{-1}\otimes \alpha_S)(i_R\otimes \psi_Y\otimes i_S)(\psi_X^{-1}\otimes \alpha_R\otimes i_S) \\
    &= (\psi_X^{-1}\otimes \psi_X)(i_{X^{\otimes m}}\otimes \alpha_R\otimes \alpha_S)
\end{align*}
This implies, by Lemma~\ref{essential2}, that $\alpha_R\otimes \alpha_S=i_{R\otimes S}$. A similar computation using Equations~\eqref{eq.2.2} and~\eqref{eq.2.3} yields $\alpha_S\otimes\alpha_R=i_{S\otimes R}$.

Now define two new $\psi$-maps by
\begin{align*}
    \tpsi_X &= \psi_X  \\
    \tpsi_Y &= \psi_Y (i_S\otimes \alpha_R)
\end{align*}

We now claim that $X$ and $Y$ are compatible shift equivalent via $(\tpsi_X,\tpsi_Y,\varphi_R,\varphi_S)$. Checking the compatible shift equivalence condition:
\begin{align*}
    (i_R\otimes\tpsi_Y)(\tpsi_X^{-1}\otimes i_R) &= (i_R\otimes \psi_Y)(i_{R\otimes S}\otimes \alpha_R)(\psi_X^{-1}\otimes i_R) \\
    &= (i_R\otimes \psi_Y)(\psi_X^{-1}\otimes \alpha_R) = \varphi_R^{(m)}.
\end{align*}

Since $\alpha_R\otimes \alpha_S=i_{R\otimes S}$, we have $\alpha_R^{-1}\otimes 1_S=i_R\otimes\alpha_S$. Therefore, 
\begin{align*}
     (i_S\otimes\tpsi_X)(\tpsi_Y^{-1}\otimes i_S) &= 
     (i_S\otimes \psi_X)(i_S\otimes \alpha_R^{-1}\otimes i_S)(\psi_Y^{-1}\otimes i_S) \\
     &= 
     (i_S\otimes \psi_X)(i_S\otimes i_R\otimes \alpha_S)(\psi_Y^{-1}\otimes i_S) \\
     &= (i_S\otimes\psi_X)(\psi_Y^{-1}\otimes \alpha_S)=\varphi_S^{(m)}.
\end{align*}
This proves that $X$ and $Y$ are compatible shift equivalent. 
\end{proof}

Let us denote
\begin{align*}
    \rho_X &= (\psi_X\otimes i_X)(i_R\otimes \varphi_S)(\varphi_R\otimes i_S)(i_X\otimes \psi_X^{-1}): X^{\otimes(m+1)}\to X^{\otimes (m+1)} \\
    \rho_Y &= (\psi_Y\otimes i_Y)(i_S\otimes \varphi_R)(\varphi_S\otimes i_R)(i_Y\otimes \psi_Y^{-1}): Y^{\otimes(m+1)}\to Y^{\otimes (m+1)}
\end{align*}
We note that aligned shift equivalence is the same as saying both $\rho_X$ and $\rho_Y$ are identity maps. 

We first prove two lemmas for compatible shift equivalences. 

\begin{lemma}\label{lm.rho.interchange} Compatible shift equivalence with lag $m$ implies that
\begin{align*}
    (i_S\otimes \rho_X)\varphi_S^{(m+1)} &= \varphi_S^{(m+1)}(\rho_Y\otimes i_S) \\
    (i_R\otimes \rho_Y)\varphi_R^{(m+1)} &= \varphi_R^{(m+1)}(\rho_X\otimes i_R)
\end{align*}
\end{lemma}

\begin{proof} We use the fact that
\[\varphi_S^{(m+1)}=(\varphi_S^{(m)}\otimes i_X)(i_{Y^{\otimes m}}\otimes \varphi_S) = (\varphi_S\otimes i_{X^{\otimes m}})(i_Y\otimes \varphi_S^{(m)}).\]
By compatible shift equivalence with lag $m$, 
\[\varphi_S^{(m)} = (i_S\otimes\psi_X)(\psi_Y^{-1}\otimes i_S).\]
We have,
\begin{align*}
    &(i_S\otimes \rho_X)\varphi_S^{(m+1)} \\
    =& (i_S\otimes \psi_X\otimes i_X)(i_{S\otimes R}\otimes \varphi_S)(i_S\otimes\varphi_R\otimes i_S)(i_{S\otimes X}\otimes\psi_X^{-1})(\varphi_S\otimes i_{X^{\otimes m}})(i_Y\otimes\varphi_S^{(m)}) \\
    =& (i_S\otimes \psi_X\otimes i_X)(i_{S\otimes R}\otimes \varphi_S)(i_S\otimes\varphi_R\otimes i_S)(\varphi_S\otimes\psi_X^{-1})\left[i_Y\otimes \left((i_S\otimes\psi_X)(\psi_Y^{-1}\otimes i_S)\right)\right] \\
    =& (i_S\otimes\psi_X\otimes i_X)(i_{S\otimes R}\otimes \varphi_S)(i_S\otimes\varphi_R\otimes i_S)(\varphi_S\otimes i_{R\otimes S})(i_Y\otimes \psi_Y^{-1}\otimes i_S) \\
    =& (\varphi_S^{(m)}\otimes i_X)(\psi_Y\otimes i_{S\otimes X})(i_{S\otimes R}\otimes \varphi_S)(i_S\otimes\varphi_R\otimes i_S)(\varphi_S\otimes i_{R\otimes S})(i_Y\otimes \psi_Y^{-1}\otimes i_S) \\
    =& (\varphi_S^{(m)}\otimes i_X)(i_{Y^{\otimes m}}\otimes \varphi_S)\left[\left((\psi_Y\otimes i_Y)(i_S\otimes\varphi_R)(\varphi_S\otimes i_R)(i_Y\otimes\psi_Y^{-1})\right)\otimes i_S\right] \\
    =& \varphi_S^{(m+1)}(\rho_Y\otimes i_S).
\end{align*}
A similar computation can be made for $(i_R\otimes \rho_Y)\varphi_R^{(m+1)}$ to show the other equality.
\end{proof}

\begin{lemma}\label{lm.comp.se.2} Compatible shift equivalence with lag $m$ implies 
\begin{align*}
    i_{X^{\otimes m}}\otimes \rho_X &= \rho_X\otimes i_{X^{\otimes m}} \\
    i_{Y^{\otimes m}}\otimes \rho_Y &= \rho_Y\otimes i_{Y^{\otimes m}} 
\end{align*}
\end{lemma}

\begin{proof} From the compatible shift equivalence, we have
\begin{align*}
    \varphi_S^{(m+1)}=(\varphi_S^{(m)}\otimes i_X)(i_{Y^{\otimes m}}\otimes \varphi_S) &= (\varphi_S\otimes i_{X^{\otimes m}})(i_Y\otimes \varphi_S^{(m)}) \\
    \varphi_R^{(m+1)}=(\varphi_R^{(m)}\otimes i_Y)(i_{X^{\otimes m}}\otimes \varphi_R) &= (\varphi_R\otimes i_{Y^{\otimes m}})(i_X\otimes \varphi_R^{(m)}).
\end{align*}
Apply the compatible conditions on $\varphi_R^{(m)}$ and $\varphi_S^{(m)}$ to the equations above, we obtain the following equalities:
\begin{align}
    \varphi_S^{(m+1)} &= (i_S\otimes\psi_X\otimes i_X)(\psi_Y^{-1}\otimes\varphi_S) \label{eq.S1} \\
    &= (\varphi_S\otimes\psi_X)(i_Y\otimes\psi_Y^{-1}\otimes i_S) \label{eq.S2}  \\
    \varphi_R^{(m+1)} &= (i_R\otimes\psi_Y\otimes i_Y)(\psi_X^{-1}\otimes\varphi_R) \label{eq.R1} \\
    &= (\varphi_R\otimes\psi_Y)(i_X\otimes\psi_X^{-1}\otimes i_R) \label{eq.R2} 
\end{align}
Now consider \[\beta_X=(\psi_X\otimes i_{X^{\otimes(m+1)}})(i_R\otimes \varphi_S^{(m+1)})(\varphi_R^{(m+1)}\otimes i_S)(i_{X^{\otimes(m+1)}}\otimes\psi_X^{-1}).\]
On one hand, apply Equation~\eqref{eq.S1} and \eqref{eq.R1}, we have,
\begin{align*}
    \beta_X &= (\psi_X\otimes i_{X^{\otimes(m+1)}})(i_{R\otimes S}\otimes \psi_X\otimes i_X)(i_R\otimes \psi_Y^{-1}\otimes \varphi_S)(i_R\otimes \psi_Y\otimes i_{Y\otimes S})(\psi_X^{-1}\otimes\varphi_R\otimes i_S)(i_{X^{\otimes(m+1)}}\otimes \psi_X^{-1}) \\
    &=i_{X^{\otimes m}} \otimes \left[(\psi_X\otimes i_X)(i_R\otimes\varphi_S)(\varphi_R\otimes i_S)(i_X\otimes\psi_X^{-1})\right] = i_{X^{\otimes m}}\otimes \rho_X
\end{align*}
On the other hand, apply Equation~\eqref{eq.S2} and \eqref{eq.R2}, we have,
\begin{align*}
    \beta_X &= (\psi_X\otimes i_{X^{\otimes(m+1)}})(i_R\otimes\varphi_S\otimes \psi_X)(i_R\otimes i_Y\otimes \psi_Y^{-1}\otimes i_S)(\varphi_R\otimes \psi_Y\otimes i_S)(i_X\otimes \psi_X^{-1}\otimes i_{R\otimes S})(i_{X^{\otimes(m+1)}}\otimes \psi_X^{-1}) \\
    &= \left[(\psi_X\otimes i_X)(i_R\otimes\varphi_S)(\varphi_R\otimes i_S)(i_X\otimes\psi_X^{-1})\right]\otimes i_{X^{\otimes m}} = \rho_X \otimes i_{X^{\otimes m}}
\end{align*}
Therefore, $i_{X^{\otimes m}}\otimes \rho_X=\rho_X \otimes i_{X^{\otimes m}}$. The other equality follows similarly. 
\end{proof}

\begin{proposition}\label{prop.c.to.a} Compatible shift equivalence with lag $m$ implies aligned shift equivalence with lag $2m$. 
\end{proposition}

\begin{proof} Define $R_1=R\otimes_B S\otimes_A R$, $S_1=S$. Here, $R_1$ is a $A$-$B$-correspondence, and $S$ is a $B$-$A$-correspondence. Define the following maps:
\begin{align*}
    \tpsi_{X} &= (\psi_X\otimes \psi_X): \quad R_1\otimes S_1 \to X^{\otimes 2m} \\
    \tpsi_{Y} &= (\psi_Y\otimes \psi_Y): \quad S_1\otimes R_1\to Y^{\otimes 2m} \\
    \tvarphi_{R} &= (\psi_X^{-1}\otimes \varphi_R)(\rho_X^{-1}\otimes i_R)(i_X\otimes\psi_X\otimes i_R): \quad X\otimes R_1\to R_1\otimes Y \\
    \tvarphi_{S} &= \varphi_{S}:\quad Y\otimes S_1\to S_1 \otimes X
\end{align*}
We claim that $X$ and $Y$ are aligned shift equivalent (with lag $2m$) via $(\tpsi_X, \tpsi_Y, \tvarphi_R, \tvarphi_S)$. 

We first compute
\[\trho_X=(\tpsi_X\otimes i_X)(i_{R_1}\otimes \tvarphi_S)(\tvarphi_R\otimes i_{S_1})(i_X\otimes \tpsi_X^{-1}).\]
From the definition,
\begin{align*}
    \trho_X &= ({\color{black}\psi_X}\otimes \psi_X\otimes i_X)(i_{R\otimes S\otimes R}\otimes \varphi_S)({\color{black}\psi_X^{-1}}\otimes \varphi_R\otimes i_S)(\rho_X^{-1}\otimes i_{R\otimes S})(i_X\otimes {\color{black}\psi_X}\otimes i_{R\otimes S})(i_X\otimes {\color{black}\psi_X^{-1}}\otimes \psi_X^{-1}) \\
    &= \left[i_{X^{\otimes m}}\otimes \left((\psi_X\otimes i_X)(i_R\otimes \varphi_S)(\varphi_R\otimes i_S)(i_X\otimes \psi_X^{-1})\right)\right] (\rho_X^{-1}\otimes i_{X^{\otimes m}}) \\
    &= (i_{X^{\otimes m}}\otimes\rho_X)(\rho_X^{-1}\otimes i_{X^{\otimes m}}).
\end{align*}
By Lemma~\ref{lm.comp.se.2}, $i_{X^{\otimes m}}\otimes \rho_X = \rho_X\otimes i_{X^{\otimes m}}$, and therefore, $\trho_X$ is identity map. 

We now compute 
\begin{align*}
    \trho_Y &=(\tpsi_Y\otimes i_Y)(i_{S_1}\otimes \tvarphi_R)(\tvarphi_S\otimes i_{R_1})(i_Y\otimes \tpsi_Y^{-1}) \\
    &= (\tpsi_Y\otimes i_Y)(i_S\otimes\psi_X^{-1}\otimes\varphi_R)(i_S\otimes\rho_X^{-1}\otimes i_R)({i_{S\otimes X}}\otimes \psi_X\otimes i_R)({\varphi_S}\otimes i_{R\otimes S\otimes R})(i_Y\otimes \tpsi_Y^{-1})\\
    &= (\tpsi_Y\otimes i_Y)(i_S\otimes\psi_X^{-1}\otimes\varphi_R)(i_S\otimes\rho_X^{-1}\otimes i_R)(\varphi_S\otimes \psi_X\otimes i_R)(i_Y\otimes \tpsi_Y^{-1}).
\end{align*}
Let us first expand the term $(i_S\otimes \rho_X^{-1})(\varphi_S\otimes\psi_X)$. 
\begin{align*}
     &(i_S\otimes\rho_X^{-1})({\varphi_S} \otimes \psi_X) \\
     =&(i_S\otimes\rho_X^{-1})({\varphi_S} \otimes i_{X^{\otimes m}})(i_{Y\otimes S}\otimes \psi_X) \\
     =& (i_S\otimes\rho_X^{-1})\varphi_S^{(m+1)}(i_Y\otimes\left(\varphi_S^{(m)}\right)^{-1})(i_{Y\otimes S}\otimes \psi_X).
\end{align*}
Apply Lemma~\ref{lm.rho.interchange} and the compatible shift equivalence condition for $\varphi_S^{(m)}$, we have,
\begin{align*}
     &(i_S\otimes\rho_X^{-1})({\varphi_S} \otimes \psi_X) \\
     =& \varphi_S^{(m+1)}(\rho_Y^{-1}\otimes i_S)(i_Y\otimes\left[(\psi_Y\otimes i_S)(i_S\otimes\psi_X^{-1})\right] (i_{Y\otimes S}\otimes \psi_X)\\
     =& (\varphi_S^{(m)}\otimes i_X)(i_{Y^{\otimes m}}\otimes \varphi_S)(\rho_Y^{-1}\otimes i_S)(i_Y\otimes \psi_Y\otimes i_S)\\
     =& ({i_S\otimes\psi_X}\otimes i_{X})(\psi_Y^{-1}\otimes i_{S\otimes  X})(i_{Y^{\otimes m}}\otimes \varphi_S)(\rho_Y^{-1}\otimes i_S)(i_Y\otimes \psi_Y\otimes i_S)
\end{align*}
Replace this term back to the earlier computation of $\trho_Y$, we have
\begin{align*}
    \trho_Y &= (\psi_Y\otimes\psi_Y\otimes i_Y)(\psi_Y^{-1}\otimes i_S\otimes\varphi_R)(i_{Y^{\otimes m}}\otimes \varphi_S\otimes i_R)(\rho_Y^{-1}\otimes i_{S\otimes R})(i_{Y^{\otimes (m+1)}}\otimes\psi_Y^{-1}) \\
    &= \left[i_{Y^{\otimes m}}\otimes \left[(\psi_Y\otimes i_Y)(i_S\otimes\varphi_R)(\varphi_S\otimes i_R)(i_Y\otimes\psi_Y^{-1})\right]\right](\rho_Y^{-1}\otimes i_{Y^{\otimes m}}) \\
    &= (i_{Y^{\otimes m}}\otimes \rho_Y)(\rho_Y^{-1}\otimes i_Y^{\otimes m}). 
\end{align*}
Again, by Lemma~\ref{lm.comp.se.2}, we have $i_{Y^{\otimes m}}\otimes \rho_Y = \rho_Y\otimes i_{Y^{\otimes m}}$ and therefore $\trho_Y$ is also the identity map. As a result, we conclude that $X$ and $Y$ are aligned shift equivalent with lag $2m$ via the maps $(\tpsi_X,\tpsi_Y,\tvarphi_R,\tvarphi_S)$. 
\end{proof}

\begin{remark} It remains unclear how aligned, balanced, and compatible shift equivalences are related to strong shift equivalence for general $\textrm{C}^*$-correspondences, even in the case of regular and full $\textrm{C}^*$-correspondences. As observed by the second named author in \cite[Theorem 6.5]{Ery22}, strong shift equivalence implies aligned shift equivalence. In particular, the Cuntz-Pimsner algebras of aligned shift equivalent $\textrm{C}^*$-correspondences are Morita equivalent. However, it remains an open question whether aligned shift equivalent $\textrm{C}^*$-correspondences are strong shift equivalent. The answer is likely in the positive, as suggested by the matrix case \cite{CDE24, BDR-SER} and recent developments on Leavitt path algebras \cite{BD2026}.
\end{remark}

\bibliographystyle{alpha}

\end{document}